\documentclass[10pt]{article}

\usepackage[centertags]{amsmath}     
\usepackage{amssymb}
\usepackage{amsthm}
\usepackage{verbatim}                
\usepackage[dvips]{graphics}
\usepackage{enumerate}

\numberwithin{equation}{section}
\newtheorem{prop}[equation]{Proposition}
\newtheorem{lem}[equation]{Lemma}
\newtheorem{dfn}[equation]{Definition}

\newtheorem{cor}[equation]{Corollary}
\newtheorem{thm}[equation]{Theorem}
\newtheorem{clm}[equation]{Claim}

\newcommand{\me}{\mathrm{e}}

\newcommand{\eps}{\varepsilon}

\newcommand{\tld}[1]{{\tilde{#1}}}
\newcommand{\vphi}{\varphi}

\newcommand{\RR}{\mathbb{R}}

\newcommand{\DIV}{\mathrm{div}}

\newcommand{\inrad}{\mathrm{inrad}}

\newcommand{\Vol}{\mathrm{Vol}}
\newcommand{\vol}{\mathrm{vol}}

\newcommand{\grad}{\mathrm{grad}}

\newcommand{\CAP}{\mathrm{cap}}
\newcommand{\Lip}{\mathrm{Lip}}

\newcommand{\ol}{\overline}

\newcommand{\la}{\langle}
\newcommand{\ra}{\rangle}
\newcommand{\bdry}{\partial}
\newcommand{\dif}[1]{{\;d #1}}
\newcommand{\dx}{\;dx}

\newcommand{\Bhf}{\frac{1}{2}B}
\newcommand{\Qhf}{\frac{1}{2}Q}

\begin{document}
\title{Local Asymmetry and the Inner Radius of Nodal Domains}
\author{Dan Mangoubi}
\date{}
\maketitle
\begin{abstract}
Let $M$ be a closed Riemannian manifold of dimension $n$. Let
$\vphi_{\lambda}$ be an eigenfunction of the Laplace--Beltrami
operator corresponding to an eigenvalue $\lambda$.
 We show  that the volume of
 $\{\vphi_\lambda >0\}\cap B$ is
 $\geq C|B|/\lambda^n$, where $B$ is any ball centered at a point of
 the nodal set.
  We apply this result to prove that each
nodal domain contains a ball of radius $\geq C/\lambda^{n}$. The
results in this paper extend previous results of F.~Nazarov,
L.~Polterovich, and M.~Sodin, and of the author.
\end{abstract}
\section{Introduction and Main Results}

\begin{figure}
\begin{center}
\scalebox{0.4}{\includegraphics{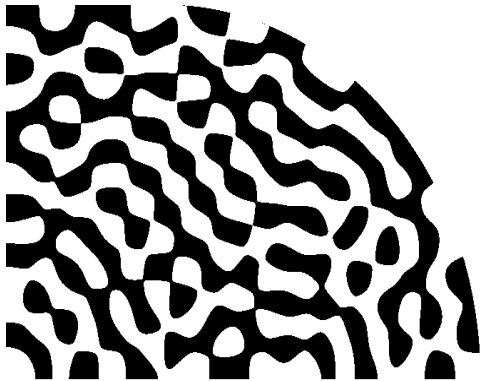}} \caption{Nodal
domains on a Quarter of a Stadium, Dirichlet boundary conditions.
Courtesy of Sven Gnutzmann}
 \label{fig:nodal}
\end{center}
\end{figure}

Let $(M, g)$ be a closed Riemannian manifold of dimension $n$, and
let $\Delta=$ \mbox{$-\DIV\circ\grad$} be the Laplace--Beltrami
operator on $M$. We consider the eigenvalue equation
\begin{equation}
\label{eigen}
\Delta \vphi_\lambda = \lambda \vphi_\lambda\ .
\end{equation}
A $\lambda$-\emph{nodal domain} on $M$ is any connected component
of the set $\{\vphi_\lambda \neq 0\}$ (see Fig.~\ref{fig:nodal}, where
the positivity set is colored in white).
In this paper we study asymptotic local geometry of nodal
domains. Let $\Omega_\lambda$ denote a $\lambda$-nodal domain on
$M$. Let $C_i$, $i=1, 2, \ldots$ denote constants which depend
only on the Riemannian metric $g$. Our first result is
\begin{thm}
\label{thm:asym}
  $$\frac{\Vol(\{\vphi_\lambda >0\}\cap B)}{\Vol(B)}\geq
  \frac{C_1}{\lambda^{(n-1)/2}}\ ,$$
  for all geodesic balls $B\subseteq M$ such that
  $\{\vphi_\lambda =0\}\cap\Bhf\neq\emptyset$. Here, $\Bhf$
  is a concentric ball of half the radius of $B$.
\end{thm}
One can think of Theorem~\ref{thm:asym} as measuring the local
asymmetry of nodal domains. Namely, it measures the volumes ratio between
the positivity and the negativity set of $\vphi_\lambda$ in $B$.
 Our motivation to prove the local asymmetry estimate
 in Theorem~\ref{thm:asym} comes from two main sources.
The first one is the following local asymmetry estimate in dimension two:
\begin{thm}[\cite{NazPolSod05}]
\label{thm:asym2}
   Let $\Sigma$ be a closed Riemannian surface.
  Then
   $$\frac{\Vol(\{\vphi_\lambda >0\}\cap B)}{\Vol(B)}\geq
   \frac{C_2}{\log\lambda\sqrt{\log\log\lambda}}\ ,$$
   for all geodesic balls $B\subseteq M$ such that
   $\{\vphi_\lambda=0\}\cap \Bhf \neq\emptyset$.
 \end{thm}
 The proof of Theorem~\ref{thm:asym2} is based on one-dimensional
 complex analysis. F.~Nazarov, L.~Polterovich and M.~Sodin suggest
 in~\cite{NazPolSod05} to explore local asymmetry in higher dimensions.
 The idea of the proof of Theorem~\ref{thm:asym} is based on
 a method of Carleman in~\cite{Car26}. Carleman finds a differential inequality
 which relates the growth of a harmonic function 
 in a two dimensional ball to its volume of positivity.
  In~\cite{NazPolSod05}, the authors indicate how
 to obtain a local asymmetry estimate for harmonic functions in
 dimensions $n\geq 3$ based on Carleman's method.
 In this paper we adapt Carleman's method to solutions of second order elliptic
 equations. As a result we can get a local asymmetry estimate also for eigenfuncions
 of the Laplace--Beltrami operator.

Our second source of motivation comes from
 our work~\cite{Man05}. In that work we gave a lower bound
 for the inner radius of nodal domains based on a growth bound
 for eigenfunctions by H.~Donnelly and C.~Fefferman and the Local Courant's Nodal Domain
 Theorem:
\begin{thm}[\cite{DonFef90, ChaMuc91}]
 \label{thm:local-courant}
   Let $M$ be a closed Riemannian manifold of dimension $n$.
   Let $\Omega_\lambda$ be a $\lambda$-nodal domain. Then
     $$\frac{\Vol(\Omega_\lambda\cap B)}{\Vol(B)}\geq
     \frac{C_3}{\lambda^{3 n^2}}\ ,$$
     for all geodesic balls $B\subseteq M$ such that
     $\Omega_\lambda\cap\Bhf\neq\emptyset$.
\end{thm}
 In the present paper Theorem~\ref{thm:asym} replaces
Theorem~\ref{thm:local-courant}. Namely, we now consider the union
of all components of the positivity set of $\vphi_\lambda$ in $B$,
while in Theorem~\ref{thm:local-courant} only one deep (i.e. which
intersects $\Bhf$) component in $B$ is considered. This lets us 
improve our estimate on the inner
radius significantly, and make the proof of our result better
suited for possible future generalizations. 
We believe that the lower bound  $C_3/\lambda^{C_4 n}$ is true also 
in Theorem~\ref{thm:local-courant}.
We prove:
\begin{thm}
\label{thm:inrad}
 $$\frac{C_5}{\lambda^{\alpha(n)}}\leq\inrad(\Omega_\lambda)\leq
 \frac{C_6}{\sqrt{\lambda}},$$
 where $\alpha(n) = \frac{1}{4} (n-1) + \frac{1}{2n}$.
\end{thm}

The proof of the upper bound and of the two dimensional case is given in~\cite{Man05}.
In this paper we assume $n\geq 3$.
\paragraph{Organization of the Paper.}
In Section~\ref{sec:prelim} we explain the principle that in small
scales compared with the wavelength~$1/\sqrt{\lambda}$ an
eigenfunction behaves like a harmonic function. In
Section~\ref{sec:elliptic} we present versions
of the Maximum Principle, the Harnack Inequality and the Mean Value Property 
for solutions of second order elliptic equations. We
give the proofs of some of these theorems in Section~\ref{sec:elliptic-prfs}. 
In Section~\ref{sec:pos-vol} we give an estimate of the volume of
positivity for solutions of the Schr\"odinger equation with small
potential in the unit ball. Our estimate will be given in terms of
the growth of the solution, and its proof is based on Carleman's method.
 In Section~\ref{sec:asym} we combine
our estimate from section~\ref{sec:pos-vol} and a growth bound by
Donnelly and Fefferman in order to prove Theorem~\ref{thm:asym}.
In section~\ref{sec:asym-inrad} we prove that the asymmetry
estimate in Theorem~\ref{thm:asym} implies the estimate on the
inner radius of a nodal domain in Theorem~\ref{thm:inrad}.

\paragraph{Acknowledgements.}
I would like to express my deep gratitude and appreciation to Leonid
Polterovich under the direction of whom this work began,
and also to Misha Sodin for his invaluable guidance and extensive discussions.

I thank F\"edor Nazarov and Misha Sodin for drawing my attention
to the asymmetry result in~\cite{NazPolSod05}, and the suggestion
to prove an estimate in higher dimensions.

I thank Dima Jakobson for his continuous encouragement and
support. I thank Pengfei Guan, Iosif Polterovich and John Toth for
discussions concerning elliptic theory, while visiting the CRM. I
thank Johannes Sj\"ostrand for
 explaining to me some estimates of the Poisson
 Kernel in a previous version of this paper.
 I thank J\"urgen Jost for discussing with me the Harnack Inequality,
 and finally I thank Sven Gnutzmann for giving me pictures of nodal
 domains he generated with his computer program.
 One of these pictures is shown
 in Fig.~\ref{fig:nodal}.

 This work was done while visiting
 the CRM, Montr\'eal, the University of Mcgill, Montr\'eal and
 the IH\'{E}S, Bures-sur-Yvette. The hospitality of these
 institutions is gratefully acknowledged.
\section{Eigenfunctions on the Wavelength Scale}
\label{sec:prelim}
In this section we explain the following principle.
\paragraph*{Principle:}
On a small scale comparable to the wavelength ($1/\sqrt{\lambda}$),
eigenfunctions behave like harmonic functions.
\\[1ex]
The above principle was extensively used in the works of H.~Donnelly,
C.~Fefferman and N.~Nadirashvili.

We may fix an atlas on $M$ for which all the transition maps
are of bounded $C^1$-norm. 
 In local coordinates, the coefficients of $g$ are
given by $g_{ij}$. The coefficients of the inverse matrix are
given by $g^{ij}$. 
In each chart we have
\begin{eqnarray}
\label{metric-bounds}
  \|g^{ij}\|_{C^1} \leq K_1 , &
  g = \det g_{ij} \leq K_2,
\end{eqnarray}
and an ellipticity bound
\begin{equation}
\label{elliptic-bound}
 g^{ij}(x) \xi_i\xi_j \geq \kappa |\xi|^2.
\end{equation}
The eigenequation~(\ref{eigen}) expressed in local coordinates is
\begin{equation}
\label{local-eigen}
 -\frac{1}{\sqrt{g}}\partial_i(g^{ij}\sqrt{g}\partial_j \vphi_\lambda)
 =\lambda\vphi_\lambda\ .
\end{equation}
We consider equation~(\ref{local-eigen}) in balls $B_r=B(0, r)$,
where $r=\sqrt{\eps_0/\lambda}$ and $\eps_0$ is a small positive
number to be chosen later.
 When we rescale it to an equation in the unit ball $B_1$, we get
\begin{equation}
\label{eigen-scaled} -\partial_i(g_r^{ij}\sqrt{g_r}\partial_j
\vphi_{\lambda, r}) = \eps_0\sqrt{g_r}\vphi_{\lambda, r} \mbox{ on
} B_1.
\end{equation}
Here, a subindex $r$ denotes a scaled function, i.e.~$f_r(x):= f(rx)$.
Since $r<1$, the bounds~(\ref{metric-bounds})
and~(\ref{elliptic-bound}) remain true
also for the rescaled metric coefficients.

Throughout this paper we let
\begin{equation}
\label{simple-notation}
 \vphi=\vphi_{\lambda, r},\quad a^{ij} = g^{ij}_r\sqrt{g_r}, \quad q=\sqrt{g_r}\ .
\end{equation}
We set
\begin{equation}
\label{L-form}
 Lu := -\partial_i (a^{ij}\partial_j u) -\eps_0 qu.
\end{equation}
Equation~(\ref{local-eigen}) takes now the form
\begin{equation}
\label{schr}
 L\vphi = 0 \mbox{ in } B_1,
\end{equation}
with the bounds
\begin{equation}
\label{coeff-bound}
  \|a^{ij}\|_{C^1(\ol{B_1})} \leq K_3, 0\leq q \leq K_4,
\end{equation}
and an ellipticity bound
\begin{equation}
\label{ellipticity}
  a^{ij}\xi_i\xi_j \geq K_5 |\xi|^2\ .
\end{equation}

If $\eps_0$ is small enough $L$ is close to be the Euclidean Laplacian
(after a linear change of coordinates) and $\vphi$ is close
to be a harmonic function.
%
%
\section{Estimates for Solutions of Elliptic Equations}
\label{sec:elliptic}
In this section we present some properties of
solutions, subsolutions and supersolutions of second order elliptic equations which
will be useful in the next sections.
$L$ is the operator given in~(\ref{L-form}) in the unit ball $B_1$.

The following theorem is a local maximum principle.
\begin{thm}[{\cite[Theorem 9.20]{GilTru83}}]
\label{thm:local-max}
Suppose $Lu \leq 0$ on $B_1$. Then
$$\sup_{B(y, r_1)} u \leq
C_1(r_1/r_2, p) \left(\frac{1}{\Vol(B(y, r_2))}
{\int_{B(y, r_2)} (u^+(x))^p\dx}\right)^{1/p},$$
for all $p>0$, whenever $0<r_1<r_2$ and $B(y, r_2) \subseteq B_1$.
\end{thm}

We will also need the weak Harnack Inequality
\begin{thm}[{\cite[Theorem 9.22]{GilTru83}}]
\label{thm:weak-harnack}
Suppose $Lu \geq -\delta$ in $B_1$, and $u\geq 0$ in $B(y, r_2)\subseteq B(0, 1)$.
Then $\exists p>0$ such that
$$ \left(\frac{1}{\Vol(B(y, r_1))}{\int_{B(y, r_1)} u^p}\right)^{1/p}\leq
C_2(r_1, r_2) \inf_{B(y, r_1)} u + C_3(r_1, r_2) |\delta|, $$
where $r_1< r_2$.
\end{thm}

We let
$$L_0 u := -\partial_i(a^{ij} \partial_j u )\ ,$$
where $a^{ij}$ are as in~(\ref{simple-notation}).
Then $L = L_0 - \eps_0 q$.
A maximum principle for $L_0$ is
\begin{thm}[{\cite[Theorem 3.7]{GilTru83}}]
\label{thm:max-L0}
Let $u$ satisfy $L_0 u \leq \delta$ on a ball $B\subseteq B_1$. Then
$$\sup_{\bdry B} u \geq \sup_B u - C_4 |\delta|\ ,$$
where $C_4$ depends only on the $C^1$-bounds and the ellipticity bounds of the
coefficients $a^{ij}$.
\end{thm}

We recall that we denote by $\vphi$ a solution
of the Schr\"odinger equation~(\ref{schr}).
As a corollary of Theorem~\ref{thm:max-L0} we obtain
the following maximum principle. Its proof is given in 
Section~\ref{sec:elliptic-prfs}.
\begin{cor}
\label{cor:max-L}
We have
  $$\sup_{\bdry B} \vphi^+ \geq 0.9 \sup_B \vphi\ ,$$
 for all balls $B\subseteq B_1$, and for all $\eps_0$ small enough.
\end{cor}

The next theorem is a Mean Value Property. Its proof is given in 
Section~\ref{sec:elliptic-prfs}.
\begin{thm}
\label{thm:mean-value-sup}
  Suppose $\vphi(0)=0$. Then
  $$\sup_{B_{r_1}} \vphi^{-} \leq C_{5}(r_1, r_2)\sup_{B_{r_2}} \vphi^{+} \ ,$$
  where $r_1<r_2\leq 1$.
\end{thm}

%
%
\section{Positivity Volume for Solutions of \\
Schr\"odinger's Equation}
\label{sec:pos-vol}

We recall that $\vphi$ is a solution of the Schr\"odinger
equation~(\ref{schr}) in the unit ball $B_1$, under the
conditions~(\ref{coeff-bound})--(\ref{ellipticity}). We estimate
the positivity volume of $\vphi$ in terms of its growth.

Let $0<r<1$. Denote by $\beta_r^+(\vphi)$ the
growth exponent of $\vphi$:
  $$\beta_r^+(\vphi) :=\log\left|\frac{\sup_{|x|\leq 1}\vphi(x)}
  {\sup_{|x|\leq r}\vphi(x)}\right|\ . $$
 Set $\la\beta_r^+\ra = \max\{\beta_r^+, 3\}$. We prove
\begin{thm}
 \label{thm:pos-vol}
 Suppose $\vphi(0)=0$ and $\eps_0$ is small enough.
Then
  $$\Vol(\{\vphi>0\}) \geq \frac{C_1(r)}{{\la\beta_r^+\ra}^{n-1}}.$$
\end{thm}

%
We start by considering the case $\vphi(0)\neq 0$.
\begin{prop}
\label{prop:vol-near>0}
  Let $|x_0| < 1$.
  Suppose $\vphi(x_0)>0$ and $\vphi(x)\leq \gamma \vphi(x_0)$ for all
  $x\in B=B(x_0, r)\subseteq B(0, 1)$.
  Then
  $$\frac{\Vol(\{\vphi>0\}\cap B)}{\Vol(B)} \geq \frac{C_2}{\gamma}\ . $$
\end{prop}
%
%
%

 \begin{proof}[Proof of Proposition~\ref{prop:vol-near>0}]
 We apply to $\vphi$ Theorem~\ref{thm:local-max}.
 \begin{equation}
 \begin{split}
 \vphi(x_0)&\leq \sup_{B(x_0, r/2)} \vphi \leq \frac{C_3}{\Vol(B)}\int_B
 \vphi^+(x)\dx
  = \frac{C_3}{\Vol(B)}\int_{B\cap\{\vphi>0\}} \vphi(x)\dx \leq \\
  &\leq \frac{C_3\gamma}{\Vol(B)}\int_{B\cap\{\vphi>0\}} \vphi(x_0)\dx
  =C_3\gamma\frac{\Vol(\{\vphi>0\}\cap B)}{\Vol(B)} \vphi(x_0)\ .
 \end{split}
 \end{equation}
 Dividing by $\vphi(x_0)$ gives us the result.
 \end{proof}

We now treat the case $\vphi(0)=0$.
\begin{proof}[Proof of Theorem~\ref{thm:pos-vol}]
%
%
%
%
%
Let $m=\lfloor\la\beta_r^+\ra\rfloor$. Decompose the annulus
\mbox{$r<|x|<1$} into $m$ annuli $r_k<|x|<r_{k+1}$, where
 $r_k = r+(1-r) k/m$  for $k=0, \ldots m$.
Define
  $$\beta_k = \log\frac{\sup_{|x|\leq r_{k+1}} \vphi(x)} {\sup_{|x|\leq r_k} \vphi(x)},
  \quad(0\leq k\leq m-1).$$
Let $S=\{k:\, \beta_k\leq 2\beta_r^+/m \}$. Observe that $\sum_k \beta_k =
\beta_r^+$. Therefore, $|S|\geq m/2$. Let $S'$ be a maximal subset of $S\setminus\{0\}$ such
that for all $k_1, k_2\in S'$ we have $|r_{k_1}-r_{k_2}|\geq
2(1-r)/m$. Notice that $|S'| \geq (m-2)/4$.

Fix $k\in S'$. By Corollary~\ref{cor:max-L}, we can find $x_k$
such that $|x_k|=r_k$ and $$\vphi(x_k)\geq 0.9\sup_{|x|\leq
r_{k}}\vphi(x) \ .$$ Consider the ball $B=B(x_k, (1-r)/m)$. For
all $x\in B$ we have $\vphi(x)\leq
\me^{2\beta_r^+/m}\vphi(x_k)/0.9$. Hence, from
Proposition~\ref{prop:vol-near>0} we know that
$$\frac{\Vol(\{\vphi>0\}\cap B)} {\Vol(B)} \geq C_4\me^{-2\beta_r^+/m}\geq C_4\me^{-2}
\geq C_5\ .$$

If we run over all $k\in S'$, we obtain the following estimate
\begin{equation*}
\begin{split}
 \Vol(\{\vphi>0\})&\geq m\Vol(\{\vphi>0\}\cap B)/4\geq
 C_5\Vol(B)m/4 \geq \\
 &\geq C_6(1-r)^n/m^{n-1}\geq
 \frac{C_6(1-r)^n}{{\la\beta_r^+\ra}^{n-1}} \ .
\end{split}
\end{equation*}
\end{proof}
\paragraph{Remark.} In the above proof if we avoid the use of the Maximum Principle,
we get a lower bound of $C(r)/\la\beta_r^+\ra^n$.
%
%

%
%


\paragraph{Different Variants of the Growth Exponent.}
We now replace $\beta_r^+$ in Theorem~\ref{thm:pos-vol}
by a more conventional growth constant:
\begin{equation}
\beta_r(\vphi) := \log\frac{\sup_{|x|\leq 1}
|\vphi(x)|}{\sup_{|x|\leq r} |\vphi(x)|}\ .
\end{equation}
We let $\la\beta_r\ra = \max\{\beta_r, 3\}$.
\begin{prop}
\label{prop:beta-gamma} Suppose $\vphi(0)=0$. Let $0<r_1<r_2<1$. Then
\begin{equation*}
 \beta_{r_1}^{+}(\vphi) \leq C_7(r_1,  r_2) \beta_{r_2}(\vphi) \ .
\end{equation*}
\end{prop}
%
%
%
%
%
\begin{proof}
The proposition amounts to proving
\begin{equation}
\label{reduce}
  \sup_{B_{r_1}} |\vphi| \leq C_8(r_1, r_2) \sup_{B_{r_2}} \vphi\ .
\end{equation}
We may assume $\sup_{B_{r_1}} |\vphi| = \sup_{B_{r_1}} \vphi^{-}$.
But then, inequality~(\ref{reduce}) is just
Theorem~\ref{thm:mean-value-sup}.
\end{proof}
An immediate consequence of Proposition~\ref{prop:beta-gamma} and
Theorem~\ref{thm:pos-vol} is
\begin{thm}
  \label{thm:pos-vol-convention}
 Suppose $\vphi(0)=0$.
Then
  $$\Vol(\{\vphi>0\}) \geq \frac{C_9(r)}{{\la\beta_r\ra}^{n-1}},$$
  for $0<r<1$ and $\eps_0$ small enough.
\end{thm}
\section{Local Asymmetry of Nodal Domains}
\label{sec:asym}
We take the positivity volume estimate in Section~\ref{sec:pos-vol},
and a growth estimate by Donnelly and Fefferman in order to prove
Theorem~\ref{thm:asym}.
%
%

\begin{proof}[Proof of Theorem~\ref{thm:asym}]
  First, we consider balls $B\subseteq M$ in scales small compared with the wavelength
  $1/\sqrt{\lambda}$,
  i.e.\ balls whose radius $r\leq \sqrt{\eps_0/\lambda}$.
  We can assume that $B$ is the Euclidean ball $B(0, r)$.
  Let $x_0$ be such that $\vphi_\lambda(x_0) = 0$ and $|x_0|<r/2$.
  We consider the eigenfunction $\vphi_\lambda$ on the ball
  $\tld{B} = B(x_0, r/2)$. We apply Theorem~\ref{thm:pos-vol-convention}
  with the function $\vphi(x) = \vphi_\lambda(rx/2)$ which is defined
  on the unit ball $B_1$. We learn that
\begin{equation}
\begin{split}
\label{vol-estm}
 \frac{\Vol(\{\vphi_\lambda>0\}\cap B)}{\Vol(B)} \geq
 \frac{\Vol(\{\vphi_\lambda>0\}\cap
 \tilde{B})}{2^n\Vol(\tilde{B})}
 &=\frac{\Vol(\{\vphi>0\} \cap B_1)}{2^n\Vol(B_1)}
  \\ &\geq \frac{C_1}{{\la\beta_{1/2}(\vphi)\ra}^{n-1}}\ .
\end{split}
\end{equation}

Next, we recall the growth estimate for eigenfunctions by
Donnelly and Fefferman:
\begin{thm}[\cite{DonFef88}]
\label{thm:df-growth}
  $\beta_{1/2}(\vphi_\lambda; \tld{B}) \leq C_2\sqrt{\lambda}$,
  where $\beta_{1/2}(\vphi_\lambda; \tld{B})$ is by definition $\beta_{1/2}(\vphi)$.
\end{thm}
Together with~(\ref{vol-estm}) we get
\begin{equation}
\label{smallball} \frac{\Vol(\{\vphi_\lambda>0\}\cap B)}{\Vol(B)}
\geq \frac{C_3}{\lambda^{(n-1)/2}}\ .
\end{equation}

We now consider large balls $B$. Let $r>\sqrt{\eps_0/\lambda}$. We
know that the inner radius of nodal domains is $<
C_4/\sqrt{\lambda}$ (see e.g.~\cite{Man05}). From this fact it
follows
\begin{lem}
  We can find a maximal set of disjoint balls $B_i=B_i(x_i, r_0)$
  contained in $B$, such that
  $r_0< \sqrt{\eps_0/\lambda}$, $\vphi_\lambda(x_i)=0$, and
  $\Vol(\cup_i B_i)/\Vol(B) \geq C_5$.
\end{lem}\qed

  The balls $B_i$ are small. Hence, by~(\ref{smallball})
  $$\Vol(\{\vphi_\lambda>0\}\cap B_i) \geq C_6\Vol(B_i)/\lambda^{(n-1)/2}\ .$$
  Summing over all balls $B_i$ gives us
  $$\Vol(\{\vphi_\lambda>0\}\cap B) \geq C_7
  \Vol(\cup_i B_i)/\lambda^{(n-1)/2}\geq  C_8\Vol(B)/\lambda^{(n-1)/2}\ ,$$
as desired.
\end{proof}
\section{Local Asymmetry implies Inner Radius Estimate}
\label{sec:asym-inrad}
In this section we prove that a local asymmetry of a domain
$\Omega\subseteq M$ implies a lower bound on its first eigenvalue.
Then, we apply this result to
a nodal domain in order to establish Theorem~\ref{thm:inrad}.
\begin{dfn}
Let $\Omega\subseteq M$ be a domain.
We say that $\Omega$ satisfies \emph{(ASym-$\alpha$)} if
 $$ \frac{\Vol(B\setminus\Omega)}{\Vol(B)}\geq \alpha .$$
 for all balls $B\subseteq M$ such that $(\Bhf\setminus\Omega)\neq\emptyset$.
\end{dfn}

We prove
\begin{thm}
\label{thm:volume-inrad}
 Let $M$ be of dimension $n\geq 3$. If $\Omega\subseteq M$
 satisfies \emph{(ASym-$\alpha$)}, then
   $$\lambda_1(\Omega)
   \geq C_1\frac{\alpha^{1-2/n}}{\inrad(\Omega)^2}\ .$$
\end{thm}

\paragraph*{Remark.} In dimension two, one can prove
that if each connected component of the complement has area $\geq
A$, then 
$\lambda_1(\Omega)\geq 
C\min(\sqrt{A}, \inrad(\Omega))/\inrad(\Omega)^3$.

\begin{proof}
We may assume that $\alpha >0$.
Let $\psi$ be the first Dirichlet eigenfunction on $\Omega$.
We extend $\psi$ by $0$ outside $\Omega$.

Let us fix a finite atlas $\{U_i, \kappa_i\}$ on $M$ as in
Section~\ref{sec:prelim}. Here $\kappa_i: U_i \to \RR^n$, are the
coordinate maps. The metric on each chart $U_i$ is comparable to
the Euclidean metric on the unit ball.
We divide $\kappa_i(U_i)$ into small non-overlapping small cubes
$Q_{ij}$ of size $h$ to be chosen later. Define the local Rayleigh
quotient by
\begin{equation}
  \label{rayleigh}
  R_{ij}(\psi) = \frac{\int_{\kappa_i^{-1}(Q_{ij})} |\nabla \psi|^2\dif(\vol)}
                {\int_{\kappa_i^{-1}(Q_{ij})} |\psi|^2\dif(\vol)} \ .
\end{equation}
 \begin{clm}
 \label{clm:local-rayleigh}
 \begin{equation}
 \label{local-rayleigh}
     R_{ij}(\psi) \leq K \lambda_1(\Omega)\ ,
 \end{equation}
    for some $i, j$, where $K$ is the number
    of charts in the atlas.
 \end{clm}
 \begin{proof}[Proof of Claim]
    Assume the contrary, i.e. for all $i, j$
    \begin{equation}
      \label{contra-rayleigh}
       \int_{\kappa_i^{-1}(Q_{ij})} |\nabla \psi|^2\,\dif(\vol) >
                K \lambda_1(\Omega) \int_{\kappa_i^{-1}(Q_{ij})}
                |\psi|^2\,\dif(\vol)\ .
    \end{equation}
 We sum up inequalities~(\ref{contra-rayleigh})
over all cubes~$Q_{ij}$.
\begin{equation}
\label{ineq:global}
\begin{split}
&\int_\Omega |\nabla\psi|^2\dif{(\vol)}  \geq
 \frac{1}{K} \sum_{i,j} \int_{\kappa_i^{-1}(Q_{ij})} |\nabla\psi|^2\dif{(\vol)} \\
 & > \lambda_1(\Omega) \sum_{i, j} \int_{\kappa_i^{-1}(Q_{ij})} |\psi|^2\dif{(\vol)} \geq
 \lambda_1(\Omega) \int_\Omega |\psi|^2\dif{(\vol)} \ .
\end{split}
\end{equation}
Hence, we obtain the following contradiction
$$\lambda_1(\Omega) = \frac{\int_{\Omega} |\nabla \psi|^2\dif(\vol)}
                {\int_{\Omega} |\psi|^2\dif(\vol)} > \lambda_1(\Omega)\ . $$
 \end{proof}

We now make a particular choice of $h$. Set $\Omega_i = \Omega\cap
U_i$, and let $r_i$ be the Euclidean inner radius of
$\kappa_i(\Omega_i)$. Let $h=8\max_i {r_i}$. We note that
 \begin{equation}
 \label{h-inrad}
  h < C_3\inrad(\Omega),
 \end{equation}
 where $C_3$ depends only on $g$ and the atlas chosen.

Take $Q=Q_{ij}$ from Claim~\ref{clm:local-rayleigh}. Let $\Qhf$ be
a concentric cube with parallel edges of size $h/2$. Since
$r_i<h/4$
\begin{equation}
  \Qhf\setminus\kappa_i(\Omega_i)\neq\emptyset \ .
\end{equation}
So, the asymmetry assumption on $\Omega$ tells us that
\begin{equation}
 \label{asym-assump}
  \frac{\Vol(Q\setminus\kappa_i(\Omega_i))}{\Vol(Q)} \geq
  C_4\alpha\ .
 \end{equation}

 Observe that the function $\psi\circ\kappa_i^{-1}$
 vanishes on the set $Q\setminus\kappa_i(\Omega_i)$.
 We now apply to $\psi\circ\kappa_i^{-1}$ the following Poincar\'e
 type inequality due to Maz'ya.
\begin{thm}[{\cite[\S10.1.2]{Maz85}}]
\label{thm:poinca-vol} Let $Q\subset \RR^n$ be a closed cube whose
edge is of length~$a$. Then,
  $$\int_Q |u|^2\dx\leq
  \frac{C_5 a^n}{\CAP_2(F, 2Q)}\int_Q |\nabla u|^2\dx$$
  for all $u\in \Lip(Q)$ and where $F=\{u=0\}$.
\end{thm}

We also recall
\begin{thm}[{\cite[\S2.2.3]{Maz85}}]
\label{thm:iso-cap}
 $\CAP_2(F, 2Q) \geq C_6\Vol(F)^{(n-2)/n}$ for $n\geq 3$.
\end{thm}

From inequality~(\ref{asym-assump}), Theorem~\ref{thm:poinca-vol},
Theorem~\ref{thm:iso-cap} and the fact that the metric $g$ is
comparable to the Euclidean metric on each chart, we immediately
obtain
\begin{equation}
\label{ineq:poincare}
 \int_{\kappa_i(Q)} |\psi|^2\,\dif(\vol) \leq C_7(\alpha) h^2
 \int_{\kappa_i(Q)} |\nabla\psi|^2\,\dif(\vol),
\end{equation}
 where $C_7(\alpha) = C_8/\alpha^{1-2/n}$.
%
Combining inequalities~(\ref{local-rayleigh})
and~(\ref{ineq:poincare}) we arrive at $\lambda_1(\Omega) \geq
C_9/(C_7(\alpha)h^2)$. To conclude, we recall
inequality~(\ref{h-inrad}).
\end{proof}

\paragraph{Application to the Inner Radius of Nodal Domains:}
\begin{proof}[Proof of Theorem~\ref{thm:inrad}]
 We notice that $\lambda_1(\Omega_\lambda) = \lambda$. This is true
 since $\vphi_\lambda$ is a Dirichlet eigenfunction for $\Omega_\lambda$
 with constant sign. We may assume $\vphi_\lambda< 0$ on $\Omega_\lambda$.
  Theorem~\ref{thm:inrad} is a consequence of
Theorem~\ref{thm:asym} and Theorem~\ref{thm:volume-inrad}, since
$B\setminus\Omega_\lambda\supseteq\{\vphi_\lambda\geq0\}$.
\end{proof}
\section{Proofs of Elliptic Estimates}
\label{sec:elliptic-prfs}
In this section we give the proofs of the elliptic estimates from
Section~\ref{sec:elliptic}.

We begin by the proof of the maximum principle.
\begin{proof}[Proof of Corollary~\ref{cor:max-L}]
If $\sup_B \vphi \leq 0$ the theorem is trivial. Otherwise, define
$w=\vphi/\sup_B \vphi$. Then $L_0 w = \eps_0 q w \leq \eps_0 q
\sup_B w \leq \eps_0 q \leq \eps_0 K_4$.

Hence, by Theorem~\ref{thm:max-L0} we know
  $$\sup_{\bdry B} w \geq \sup_B w - C_1 K_4\eps_0 \geq 1 - C_2 \eps_0. $$
  Hence, for all $\eps_0$  small enough we have
  $\sup_{\bdry B} w \geq 0.9$, from which we conclude
  $\sup_{\bdry B} \vphi \geq 0.9 \sup_B \vphi$.
\end{proof}

We now come to the proof of the Mean Value Property:
 \begin{proof}[Proof of Theorem~\ref{thm:mean-value-sup}]
 Let $M=\sup_{B_{r_2}} \vphi^{+}$.
 Observe that $L(M-\vphi) = LM-L\vphi  = -\eps_0 q M$.
 Hence,
 $$ -\eps_0 K_4 M \leq L(M-\vphi) \leq 0 \ .$$
 By Theorem~\ref{thm:weak-harnack}
 we have for some $p>0$,
 \begin{multline}
  \label{weak-harnack}
 \left(\frac{1}{\vol(B_{(r_1+r_2)/2})}
 \int_{B_{(r_1+r_2)/2}}
 (M-\vphi)^p\right)^{1/p}
  \leq \\
  \leq C_3(r_1, r_2)(M+\inf_{B_{(r_1+r_2)/2}} (-\vphi)) \leq C_3(r_1, r_2) M \ ,
  \end{multline}
  where the last inequality is true since $\vphi(0)=0$.
 By Theorem~\ref{thm:local-max}
 we know that
 \begin{equation}
   \label{local-max}
   \sup_{B_{r_1}} (M-\vphi) \leq
   C_4(r_2/r_1, p)\left(\frac{1}{\vol(B_{(r_1+r_2)/2})}
   \int_{B_{(r_1+r_2)/2}} (M-\vphi)^p\right)^{1/p}\ .
 \end{equation}
 Combining~(\ref{weak-harnack}) and~(\ref{local-max}) we obtain
 \begin{equation}
   \sup_{B_{r_1}} (M-\vphi) \leq C_5(r_1, r_2) M \ .
 \end{equation}
 Recalling the definition of $M$ we get
   $\sup_{B_{r_1}} \vphi^{-} \leq C_5(r_1, r_2) \sup_{B_{r_2}}
   \vphi^{+}$.
\end{proof}
%
%

%
\providecommand{\bysame}{\leavevmode\hbox to3em{\hrulefill}\thinspace}
\providecommand{\MR}{\relax\ifhmode\unskip\space\fi MR }
\providecommand{\MRhref}[2]{%
  \href{http://www.ams.org/mathscinet-getitem?mr=#1}{#2}
}
\providecommand{\href}[2]{#2}


\bigskip\par\noindent
Dan Mangoubi,\\
IH\'{E}S, Le Bois-Marie,\\
35, Route de Chartres,\\
F-91440 Bures-sur-Yvette,\\
France\\
\smallskip
\texttt{\small mangoubi@ihes.fr}
\end{document}